\documentclass[12pt]{article}

\usepackage{amsthm}
\usepackage{amsmath}
\usepackage{bbm}
\usepackage{amssymb} 

\newcommand{\Exp}{{\rm I\hspace{-0.8mm}E}}
\newcommand{\Prob}{{\rm I\hspace{-0.8mm}P}}

\newcommand{\iz}{{\rm \rlap Z\kern 2.2pt Z}}

\newcommand{\RL}{{\rm I\hspace{-0.8mm}R}}

\newtheorem{theorem}{Theorem}

\newtheorem{proposition}{Proposition}

\newtheorem{remark}{Remark}

\title{\bf\LARGE {Explicit formula for the supremum distribution
of a spectrally negative stable process}}
\author{Zbigniew Michna\\
Department of Mathematics and Cybernetics\\
Wroc{\l}aw University of Economics\\
Wroc{\l}aw}
\date{}

\begin{document}

\maketitle

\bibliographystyle{abbrv}

\begin{abstract}
In this article we get simple explicit formulas for $\Exp\sup_{s\leq t}X(s)$ where $X$ is a spectrally positive or negative L\'evy process with infinite variation. As a consequence we derive a generalization of the well-known formula for the supremum distribution of Wiener process that is we obtain \\
$\Prob(\sup_{s\leq t}Z_{\alpha}(s)\geq u)=\alpha\,\Prob(Z_{\alpha}(t)\geq u)$ for $u\geq 0$ where 
$Z_{\alpha}$ is a spectrally negative L\'evy process with $1<\alpha\leq 2$ which also stems from Kendall's identity for the first crossing time. Our proof uses a formula for the supremum distribution
of a spectrally positive L\'evy process which follows easily from the elementary Seals formula.

\vspace{5mm}
{\it Keywords: L\'evy process, distribution of the supremum of a stochastic process,
$\alpha$-stable L\'evy process}
\newline
\vspace{2cm}
MSC(2010): Primary 60G51; Secondary 60G52, 60G70.
\end{abstract}

\section{Introduction}
L\'evy processes appear in many theoretical  and practical fields where they serve as a basic skeleton for a description of certain phenomena. They are applied in physics, economics, finance, insurance, queueing systems and other branches of knowledge. Their features like independence and stationarity of increments or self-similarity in certain cases permit to apply them to model for instance returns of stock prices, claims to insurance companies or an inflow (outflow) to the buffer in queueing (telecommunications) systems. Moreover L\'evy processes serve as a starting point for more complicated models e.g. based on stochastic differential equations. 

We will investigate real valued L\'evy processes.
L\'evy-It\^o representation shows their stochastic construction which is the following (see e.g. Sato \cite{sa:99})
$$
X(t)=B(t)+\int_{|x|<1}x\,(N_t(dx)-tQ(dx))+\int_{|x|\geq 1}x\,N_t(dx)+at\,,
$$
where $B(t)$ is Wiener process, $N$ is a point process generated by the jumps of $X$ that is $N=\sum_{\{t:\Delta X(t)\neq 0\}}\delta_{(t,\Delta X(t))}$. 
$N$ is a random Poisson measure  on $[0,\infty)\times\{\RL\setminus 0\}$ with the mean $ds\times Q(dx)$,
where $Q(dx)$ is the so-called L\'evy measure on $\RL\setminus 0$ and $a\in\RL$.

In this note we consider spectrally one-sided L\'evy processes without Wiener component. We find expected value of 
the supremum on a finite interval for any spectrally positive or negative L\'evy process.
Then as a corollary we derive a generalization of the famous formula
$$
\Prob(\sup_{s\leq t}B(s)\geq u)=2\Prob(B(t)\geq u)
$$
where $B$ is Wiener process that is we show that
\begin{equation}\label{mfor}
\Prob(\sup_{s\leq t}Z_{\alpha}(s)\geq u)=\alpha\Prob(Z_{\alpha}(t)\geq u)
\end{equation}
where $Z_{\alpha}$ is an $\alpha$-stable L\'evy process with $1<\alpha\leq 2$ the skewness parameter $\beta=-1$ and the shift parameter $\mu=0$ (see e.g. Janicki and Weron \cite{ja:we:94} or Samorodnitsky and Taqqu \cite{sa:ta:94}).  The formula (\ref{mfor}) also stems from Kendall's identity for the first crossing time see Kendall \cite{ke:57} or e.g. Bertoin \cite{be:96a} or Borovkov and Burq \cite{bo:bu:01} and the references therein. Let us recall that the proofs of Kendall's identity are analytical (using Laplace transforms) or are using limit and combinatorial arguments or factorization identities except the proof of Borovkov and Burq \cite{bo:bu:01} which is straightforward by the change of measure technique. 
The above formula for Wiener process follows easily from the reflection principle. Here we give a straightforward  proof based on the formula from Michna \cite{mi:11} and \cite{mi:11a} which is simply derived
from the elementary Seals formula for a compound Poisson process see Seal \cite{se:74} (which is going back to Cram\'er and Prabhu see e.g. Asmussen and Albrecher \cite{as:al:00} and Prabhu \cite{pr:61}).
Regardless the theoretical importance of the above formula, the supremum distribution is the key value in many practical problems in insurance, finance and queueing systems.  The distribution of the supremum of spectrally one-sided L\'evy processes has been investigated in many papers see e.g. Albin \cite{al:93},
Avram et al. \cite{av:ky:pi:04}, Bernyk et al. \cite{be:da:pe:08}, Bertoin \cite{be:96},
Michna \cite{mi:11}, Pistorius \cite{pi:04} and many others. Explicit formulas for the supremum distribution of stochastic processes on finite intervals are known only in few cases. Most papers are concerned with an asymptotic behavior of the tail distribution
of the supremum for stochastic processes see e.g. Albin and Sunden \cite{al:su:09} and the references therein. In some articles one can find the distribution of the supremum
but in the form of Laplace transforms of the first passage times see Bertoin \cite{be:96},
Avram et al. \cite{av:ky:pi:04} and Pistorius \cite{pi:04}.

\section{Expected value of the supremum}
Let $X$ be a spectrally positive L\'evy process and $Y$ a spectrally negative L\'evy process both with infinite variation (one can regard that $Y=-X$). Let us recall that a spectrally positive L\'evy process has no negative jumps and analogically for a spectrally negative L\'evy process. Additionally we assume that their L\'evy measure $Q$ has a bounded density on every infinite interval cut off from zero and their one-dimensional distributions are absolutely continuous with respect to Lebesgue measure see Michna \cite{mi:11a}. We denote $x^+=\max(x,0)$ and $x^-=-\min(x,0)$.
\begin{proposition}\label{expsup}
\begin{eqnarray*}
\lefteqn{\Exp\sup_{s\leq t}X(s)}\\
&=&\int_0^{\infty}\Prob(X(t)>u)\,du+\int_0^t\frac{P(X(t-s)>0)}{s}\,ds
\int_{-\infty}^{0}\Prob(X(s)\leq u)\,du\,.
\end{eqnarray*}
If $\Exp Y^{-}(t)<\infty$ then
\begin{eqnarray*}
\lefteqn{\Exp\sup_{s\leq t}Y(s)}\\
&=&\int_0^{\infty}\Prob(Y(t)>u)\,du+\int_0^t\frac{P(Y(t-s)<0)}{s}\,ds
\int_{0}^{\infty}\Prob(Y(s)\geq u)\,du\,,
\end{eqnarray*}
\end{proposition}
\proof
By Michna \cite{mi:11} and \cite{mi:11a} we have
\begin{equation}\label{fmi}
\Prob(\sup_{s\leq t}X(s)>u)=\Prob(X(t)>u)+\int_0^t\frac{f(u,s)}{t-s}\,ds\int_{-\infty}^0
\Prob(X(t-s)\leq x)\,dx\,,
\end{equation}
where $f(u,s)$ is a density function of the random variable $X(s)$.
Integrating we get
\begin{eqnarray*}
\lefteqn{\Exp\sup_{s\leq t}X(s)}\\
&=& \int_0^\infty\Prob(X(t)>u)\,du+\int_0^\infty du\int_0^t\frac{f(u,s)}{t-s}\,ds\int_{-\infty}^0
\Prob(X(t-s)\leq x)\,dx\\
&=& \int_0^\infty\Prob(X(t)>u)\,du+\int_0^t\frac{\Prob(X(s)>0)}{t-s}\,ds\int_{-\infty}^0
\Prob(X(t-s)\leq x)\,dx\\
&=& \int_0^\infty\Prob(X(t)>u)\,du+\int_0^t\frac{\Prob(X(t-s)>0)}{s}\,ds\int_{-\infty}^0
\Prob(X(s)\leq u)\,du
\end{eqnarray*}
where in the last equality we substitute $s'=t-s$.

To prove the second assertion let us notice that for a fixed $t$ and $0\leq s\leq t$ we have
$X(s)\stackrel{d}{=}X(t)-X(t-s)$ in the sense of finite dimensional distributions. Thus
\begin{eqnarray*}
\sup_{s\leq t}X(s)&\stackrel{d}{=}&\sup_{s\leq t}(X(t)-X(t-s))\\
&=&X(t)-\inf_{s\leq t}X(s)
\end{eqnarray*}
where the equality in distribution is in the sense of the one-dimensional distribution. Hence
$$
\sup_{s\leq t}X(s)\stackrel{d}{=}X(t)+\sup_{s\leq t}Y(s)
$$
where $Y=-X$. So by the first formula of the proposition we obtain
\begin{eqnarray*}
\Exp\sup_{s\leq t}Y(s)&=& -\Exp X(t)+\Exp\sup_{s\leq t}X(s)\\
&=&\Exp Y(t)+\Exp\sup_{s\leq t}X(s)\\
&=&\Exp Y(t)+\int_0^{\infty}\Prob(X(t)>u)\,du+\\
&&\,\,\,\,\,\,\int_0^t\frac{P(X(t-s)>0)}{s}\,ds
\int_{-\infty}^{0}\Prob(X(s)\leq u)\,du\\
&=&\Exp Y(t)+\int_0^{\infty}\Prob(Y(t)<-u)\,du+\\
&&\,\,\,\,\,\,\int_0^t\frac{P(Y(t-s)<0)}{s}\,ds
\int_{-\infty}^{0}\Prob(Y(s)\geq -u)\,du\\
&=&\int_{0}^{\infty}\Prob(Y(t)>u)\,du+\\
&&\,\,\,\,\,\,\int_0^t\frac{P(Y(t-s)<0)}{s}\,ds
\int_{0}^{\infty}\Prob(Y(s)\geq u)\,du
\end{eqnarray*}
which finishes the proof.
\begin{remark}
One can write the first formula of Prop. \ref{expsup} as
$$
\Exp\sup_{s\leq t}X(s)
=\Exp X^+(t)+\int_0^t\frac{P(X(t-s)>0)}{s}\,\Exp X^-(s)\,ds
$$
and the second formula as
$$
\Exp\sup_{s\leq t}Y(s)
=\Exp Y^+(t)+\int_0^t\frac{P(Y(t-s)<0)}{s}\,\Exp Y^+(s)\,ds\,.
$$
\end{remark}
\begin{remark}
The formulas of Prop. \ref{expsup} are valid for Wiener process as well because the formula
(\ref{fmi}) is true for Wiener process see Michna \cite{mi:11}.
\end{remark}
\section{The supremum distribution of a spectrally negative stable L\'evy process}
Now let us consider a spectrally negative $\alpha$-stable L\'evy process $Z_\alpha$ with $1<\alpha\leq 2$ (see e.g. Janicki and Weron \cite {ja:we:94} or Samorodnitsky and Taqqu \cite{sa:ta:94}). A simple proof
of  a generalization of the famous formula for the supremum distribution of Wiener process we get by
Prop. \ref{expsup} which in fact follows from the formula for the supremum distribution of a spectrally positive L\'evy process (see Michna \cite{mi:11a}).
\begin{theorem}
Let $u\geq 0$ and $Z_\alpha$ be an $\alpha$-stable L\'evy process with the skewness parameter $\beta=-1$. Then
$$
\Prob(\sup_{s\leq t}Z_{\alpha}(s)\geq u)=\alpha\Prob(Z_{\alpha}(t)\geq u)\,.
$$ 
\end{theorem} 
\proof
Let us note that $\Prob(Z_{\alpha}(s)>0)=1/\alpha$ and 
$Z_\alpha(as)\stackrel{d}{=}a^{1/\alpha}Z(s)$ for $a>0$ in the sense of finite dimensional distributions (the self-similarity property) see e.g. Samorodnitsky and Taqqu \cite{sa:ta:94}. Thus by Prop.
\ref{expsup} we have
\begin{eqnarray*}
\lefteqn{\Exp\sup_{s\leq t}Z_{\alpha}(s)}\\
&=&\int_0^{\infty}\Prob(Z_{\alpha}(t)>u)\,du+\int_0^t\frac{P(Z_{\alpha}(t-s)<0)}{s}\,ds
\int_{0}^{\infty}\Prob(Z_{\alpha}(s)\geq u)\,du\\
&=&\int_0^{\infty}\Prob(Z_{\alpha}(t)>u)\,du+\frac{\alpha-1}{\alpha}\int_0^t\frac{ds}{s}
\int_{0}^{\infty}\Prob(Z_{\alpha}(t)\geq ut^{1/\alpha}/s^{1/\alpha})\,du\\
&=&\int_0^{\infty}\Prob(Z_{\alpha}(t)>u)\,du+\frac{\alpha-1}{\alpha}\int_0^t\frac{s^{1/\alpha-1}}{t^{1/\alpha}}\,ds
\int_{0}^{\infty}\Prob(Z_{\alpha}(t)\geq u)\,du\\
&=&\int_0^{\infty}\Prob(Z_{\alpha}(t)>u)\,du+(\alpha-1)\int_{0}^{\infty}\Prob(Z_{\alpha}(t)\geq u)\,du\\
&=&\alpha\int_{0}^{\infty}\Prob(Z_{\alpha}(t)\geq u)\,du\,,
\end{eqnarray*}
where in the second equality we use the self-similarity property and in the third equality we substitute $u'=ut^{1/\alpha}/s^{1/\alpha}$.
By eq. 1.2 of Albin \cite{al:93} or Lemma 3 of Furrer et al. \cite{fu:mi:we:97} the following upper bound we can state
\begin{equation}\label{fur}
\Prob(\sup_{s\leq t}Z_{\alpha}(s)\geq u)\leq \alpha\Prob(Z_{\alpha}(t)\geq u)\,.
\end{equation}
Since the last calculations give
$$
\int_0^\infty\Prob(\sup_{s\leq t}Z_{\alpha}(s)\geq u)\,du=\int_0^\infty\alpha\Prob(Z_{\alpha}(t)\geq u)\,du
$$
thus by eq. (\ref{fur}) and the continuity and monotonicity with respect to $u$ we obtain the assertion of the theorem.
\begin{remark}
The asymptotic behavior of $\Prob(Z_{\alpha}(t)\geq u)$ for $u\rightarrow\infty$ can be found in Samorodnitsky and Taqqu 
\cite{sa:ta:94}, eq. 1.2.11 which is not regularly varying but Weibullian. 
\end{remark}
\begin{remark}
In Albin \cite{al:93} the exact asymptotic for $\Prob(\sup_{s\leq t}Z_{\alpha}(s)\geq u)$ as $u\rightarrow\infty$ has been derived in the form $C_\alpha\Prob(Z_{\alpha}(t)\geq u)$ where it was shown that $C_\alpha>1$ for $1<\alpha< 2$. Thus
we get that $C_\alpha=\alpha$ if $1<\alpha\leq 2$.
\end{remark}
\begin{remark}
The supremum distribution of a spectrally positive L\'evy process is qualitatively different than in the spectrally negative case because in the first case the supremum is attained by a jump see e.g. Bernyk et al. \cite{be:da:pe:08}. 
\end{remark}

\bibliographystyle{plainnat}

\end{document}